\documentclass[a4paper,10pt]{amsart}
\usepackage[german,english]{babel}
\usepackage[utf8x]{inputenc}
\usepackage{ucs}
\usepackage{amsmath}
\usepackage{amsfonts}
\usepackage{amssymb}
\usepackage{amsthm}
\usepackage{fontenc}
\usepackage{graphicx}
\usepackage{url}
\usepackage{color}
\usepackage[T1]{fontenc}
\usepackage{a4wide}
\usepackage{float}

\title[Parameter optimization for low-rank matrix recovery in hyperspectral imaging]{
Parameter optimization for low-rank matrix recovery in hyperspectral imaging}

\author{M. Wolfmayr}

\address[M. Wolfmayr]{Jamk University of Applied Sciences, Finland, and University of Jyv\"{a}skyl\"{a}, Finland}

\email{monika.wolfmayr@jamk.fi}

\begin{document}

\begin{abstract}
An approach to parameter optimization for the low-rank matrix recovery method in hyperspectral imaging is discussed. We formulate an optimization problem with respect to the initial parameters of the low-rank matrix recovery method. The performance for different parameter settings is compared in terms of computational times and memory. The results are evaluated by computing the peak signal-to-noise ratio as a quantitative measure.
The potential improvement of the performance of the noise reduction method is discussed when optimizing the choice of the initial values.
The optimization method is tested on standard and openly available hyperspectral data sets including Indian Pines, Pavia Centre, and Pavia University.
\end{abstract}

\maketitle


\section{Introduction}
\label{sec:introduction}

\begin{figure}[h]
    \centering
    \includegraphics[scale=0.5]{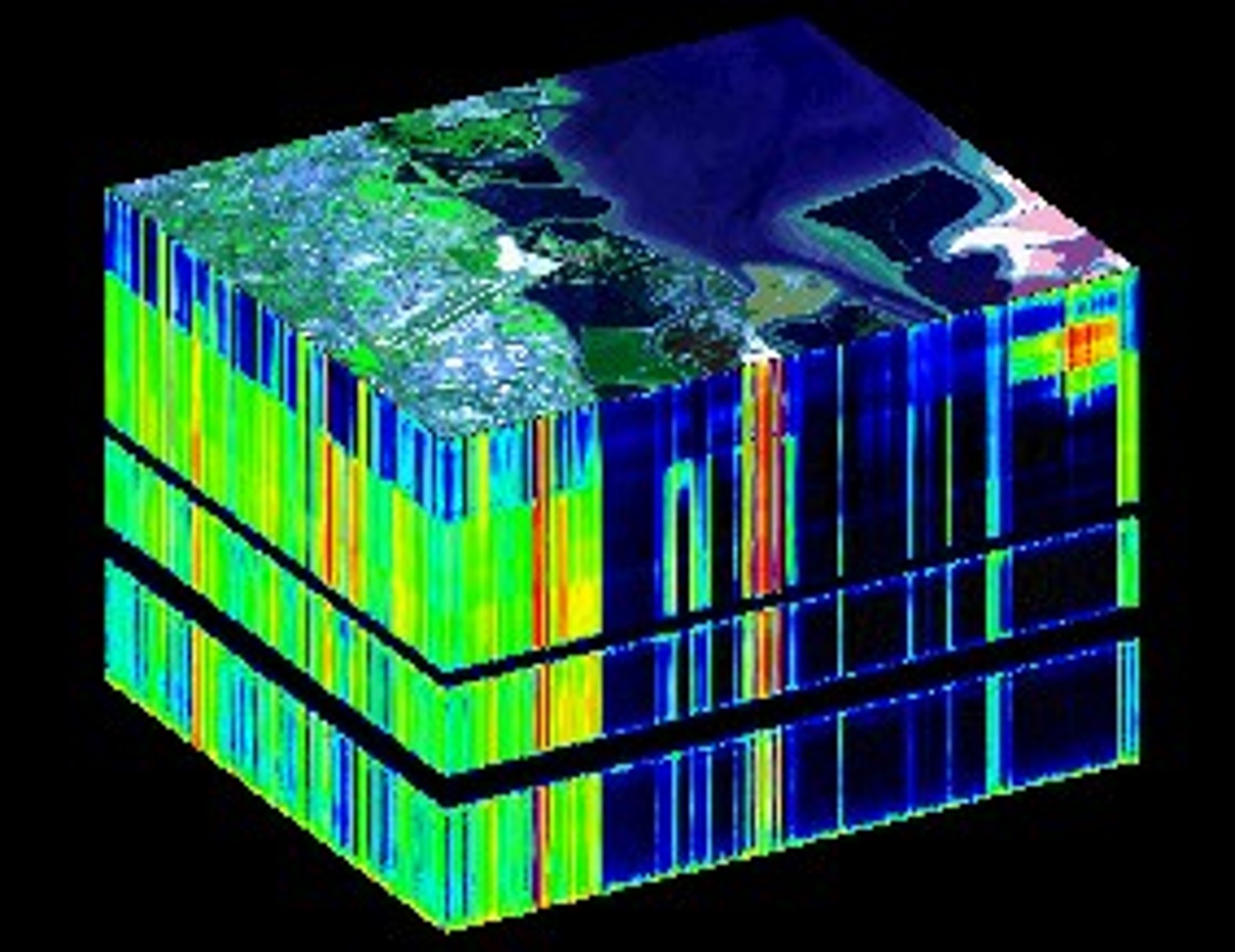}
    \caption{An example of a hyperspectral data cube (Courtesy NASA/JPL-Caltech).}
    \label{fig:datacube}
    \end{figure}

    In hyperspectral imaging (HSI), spectral signatures of objects are recorded for each image pixel. 
    The spectral signature of an object is the reflectance variation or function with respect to the wavelength. It is important to characterize materials and their properties.
    In HSI, different types of noises appear due to environmental or instrumental influences. The noises include Gaussian noise \cite{Mandelbrot}, impulse noise \cite{Majumdar}, dead pixels or lines \cite{Shen}, and stripes \cite{Rogass}. 
    This has been recently discussed for instance in \cite{Sun} and in the overview articles \cite{Ghamisi} and \cite{Rasti}.
    HSI combines spatial and spectral information in a hyperspectral data cube as displayed in Figure \ref{fig:datacube}.
    Naturally, the amount of generated data is huge, and an efficient and reliable approach to noise reduction takes advantage of the internal dependencies between the wavebands.
    HSI precision and reliability are essential for many applications including digitalization and robotization in Earth and space exploration. The applications include agriculture monitoring \cite{Singh}, atmospheric science \cite{Calin}, geology \cite{Peyghambari}, and space exploration \cite{Qian}.
    Efficient and reliable noise reduction techniques are essential for image processing in practice regarding real-time decision-making and automation.
    In \cite{Wolfmayr}, various noise reduction techniques have been compared for hyperspectral image data in asteroid imaging including also low-rank matrix recovery (LRMR).
    
    The LRMR is a low-rank modeling approach \cite{Zhang} and has been discussed among other advanced image processing methods in more detail in \cite{Ghamisi} and \cite{Rasti}.
    We use here the LRMR together with the Go Decomposition (GoDec) numerical optimization algorithm as presented in \cite{Zhou} in the inner iteration steps of the approach for solving the restoration model.
    Parameter optimization in advanced image processing can provide important indirect information for control and real-time decision-making.
    The main idea here is to optimize the choice of the algorithm's initial values in order to improve the performance of the noise reduction.
    We analyze the parameter choices including applying nonlinear optimization methods as derived in \cite{Lagarias}.
    The open available hyperspectral data sets Indian Pines \cite{IndianPines}, Pavia Centre \cite{Pavia}, and Pavia University \cite{Pavia} are used for the computational tests.

    This work is part of \textit{coADDVA - ADDing VAlue by Computing in Manufacturing} project funded by the Regional Council of Central Finland/Council of Tampere Region and European Regional Development Fund. 
    It supports the project's goals to improve the efficiency of robotics by developing optimal control methods leading to flexible imaging and automation in image processing.

    The article is organized as follows: In Section \ref{sec:methods}, the methods are discussed including the LRMR and the optimization applied to its initial parameter setting. The used data sets are presented. Section \ref{sec:results} presents the computational results. 
    Section \ref{sec:discussion} discusses the results according to existing literature and 
    possible future work. 
    Finally, Section \ref{sec:conclusions} presents the main contributions of this work.


\section{Methods and materials} 
\label{sec:methods}

\subsection{Methods} 

The first LRMR model was proposed in \cite{Wright}  as a robust principle component analysis approach. It has been further developed in \cite{Zhang} for hyperspectral image restoration by combining nonlinear optimization in the inner iteration loop in order to solve the restoration model. In this work, we apply optimization with respect to the initial parameters in terms of an outer iteration loop.

In the following, we present the LRMR model.
Given the real matrix $D$ of size $m \times n$ containing the observed data and assuming corruption by the sparse error matrix $S$ and a random Gaussian noise modeled by the matrix $N$ 
the goal is to recover the low-rank matrix $L$ 
with $D = L + S + N$, all real number matrices of the same size $m \times n$.
    The minimization problem 
    \begin{equation}
    \label{eqn:lrmr}
    \min_{L,S} \| D-L-S \|_F^2 \text{ s.t. } \text{rank}(L) \leq r, \text{card}(S) \leq p
    \end{equation}
    is solved with $r$ denoting the upper bound for the rank of $L$ and $p$ for the cardinality of $S$ which is related to the estimation of noise corruption. The norm $\| \cdot \|_F$ denotes the Frobenius-norm. The formulation \eqref{eqn:lrmr} of LRMR can be found in \cite{Zhou, Zhang}.
    Redundancies between the wavebands yield the low-rank property required for LRMR. LRMR modeling is then applied together with the GoDec algorithm \cite{Zhou} in order to solve the subproblems exploiting the low-rank property of HSI. The subproblems are created by taking subcubes centered at a pixel in the spatial dimension. Thus, if the whole data cube is of size $m \times n \times w$, where $w$ denotes the number of spectral reflectance bands, then the subcubes are of size $b \times b \times w$ with $b < m$ and $b < n$. We define the entries of the subcubes spectral band by spectral band. The entries of each subcube band- by-band are then organized in lexicographical order to obtain two-dimensional matrices of size $b^2 \times w$. The subcubes are then processed iteratively providing local image restoration. We denote further the iteration stepsize by $s$.
    More details on the specific LRMR model can be found in \cite{Zhang} and on GoDec in \cite{Zhou}.

    The main focus of this work is the detailed investigation of the LRMR method with respect to the starting values for its main variable parameters including rank $r$ and blocksize of the subcubes $b$, 
    estimation parameter for the percentage of noise corruption $p$ and stepsize $s$ of the inner iteration. In the following, we denote by $r$, $p$, $b$, and $s$ the initial setting.
    We apply nonlinear 
    optimization in order to determine the best parameter values for the parameter setting with respect to the peak signal-to-noise ratio (PSNR).
    We choose the Nelder-Mead simplex algorithm as presented in \cite{Lagarias}.
    for the nonlinear optimization.
    The PSNR is computed by
    \begin{equation}
    \label{eqn:psnr}
    PSNR = 10 * \log_{10} \frac{\max(C)^2*m*n*w}{\|C-\tilde{C}\|^2},
    \end{equation}
    where $C$ and $\tilde{C}$ denote the original and denoised data cube, respectively. The sizes of the spatial and spectral dimensions are denoted by $m, n$, and $w$.
    The performance of LRMR is analyzed in terms of computational efficiency and memory with regard to different parameter choices. 
    The computational time taken by the algorithm in Matlab and Python is compared in order to study the
    performance of LRMR and differences in the implementation between the two programming languages.
    The goal was to investigate possible difficulties in parameter optimization.

    \subsection{Materials}

    We present the data sets used in this work.

    \subsubsection{Indian Pines data set}

    The Indian Pines data set covers the Indian Pines test site which is located in North-western Indiana. It shows mainly agriculture and forest. The data set consists of 224 spectral reflectance bands and 145 $\times$ 145 pixels. The wavelength range lies betwen 0.4 and 2.5 $\times 10^{-6}$ meters.
    The Indian Pines data set is openly available in \cite{IndianPines}.
    \begin{figure}[h]
    \centering
        \includegraphics[scale=0.4]{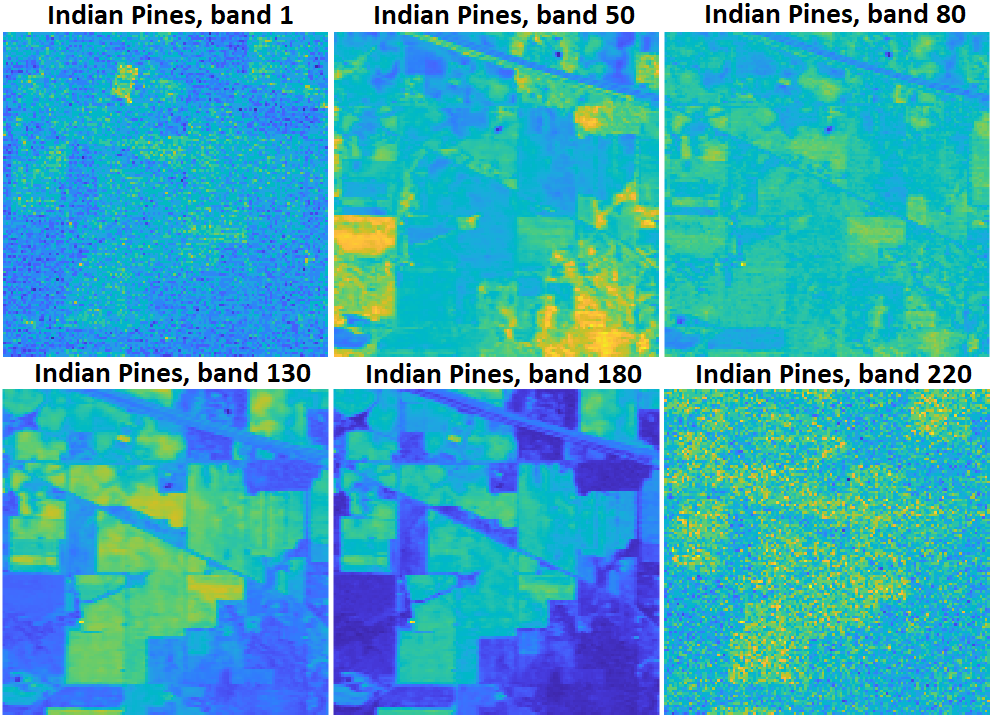}
        \caption{The Indian Pines data set displayed for six different spectral reflectance bands 1, 50, 80, 130, 180, and 220.}
        \label{fig:indian}
    \end{figure}
    Figure \ref{fig:indian} shows the scene for six different spectral reflectance bands 1, 50, 80, 130, 180, and 220. Different bands show different layers of material's spectral signature.

    \subsubsection{Pavia Centre data set}

    The data set contains scenes from the center of Pavia in northern Italy. It consists of 102 spectral reflectance bands and 1096 $\times$ 1096 pixels. The set is openly available in \cite{Pavia}.
    \begin{figure}[h]
    \centering
    \includegraphics[scale=0.35]{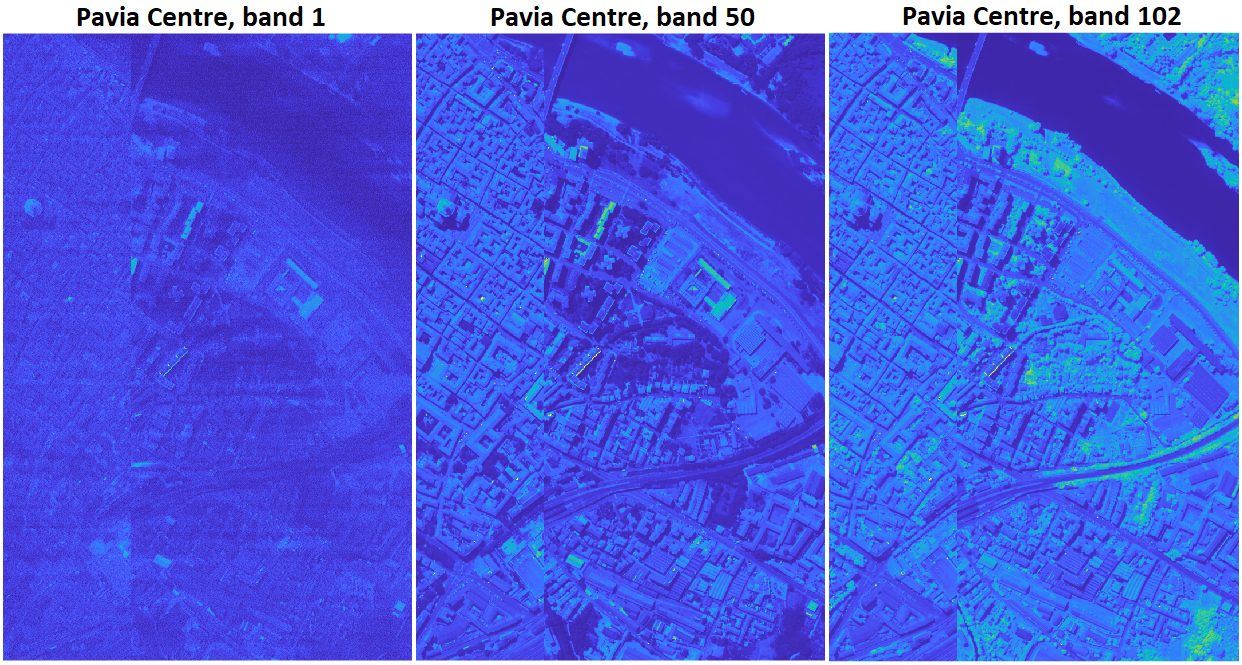}
    \caption{The Pavia Centre data set displayed for spectral reflectance bands 1, 50, and 102.}
    \label{fig:pavia}
    \end{figure}
    Figure \ref{fig:pavia} shows the scene for three different spectral reflectance bands 1, 50 and 102.

    \subsubsection{Pavia University data set}

    The data set contains scenes from the university in Pavia in northern Italy. It consists of 103 spectral reflectance bands and 610 $\times$ 610 pixels. The set is openly available in \cite{Pavia}.
    \begin{figure}[h]
    \centering
        \includegraphics[scale=0.35]{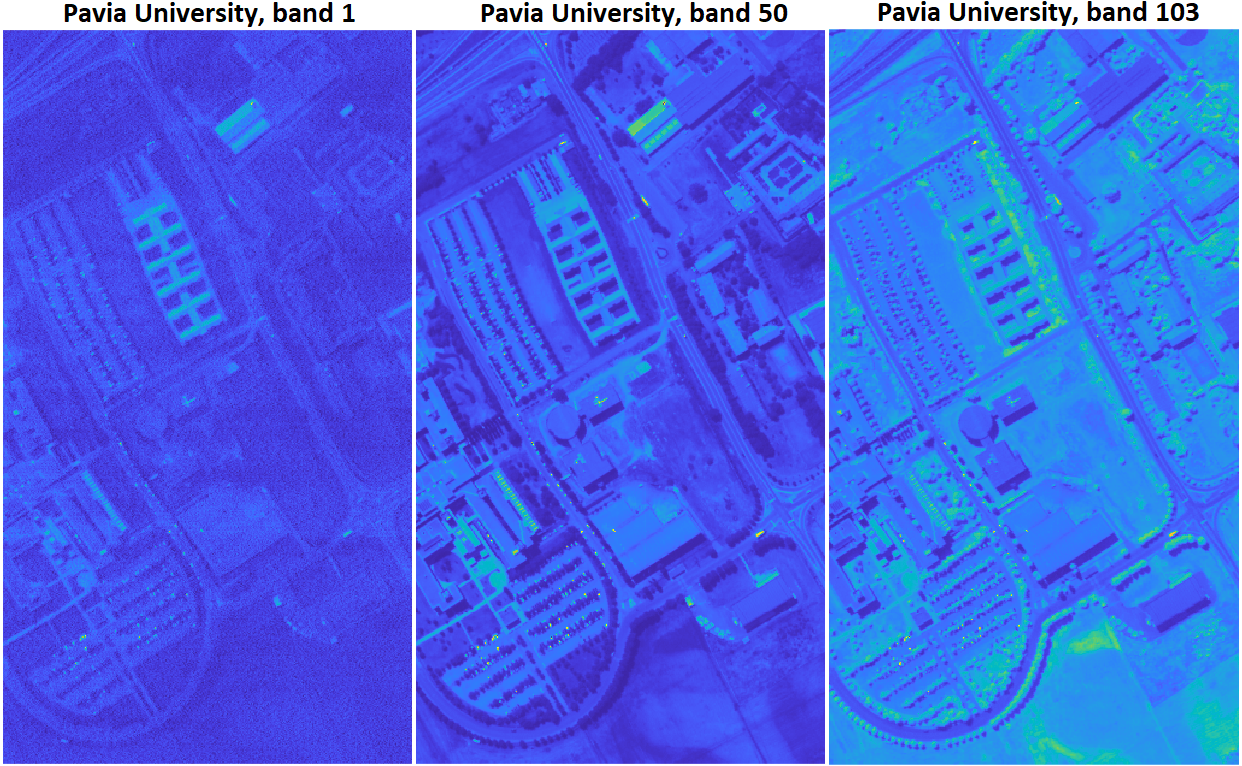}
       \caption{The Pavia University data set displayed for spectral reflectance bands 1, 50, and 103.}
       \label{fig:paviaU}
    \end{figure}
    Figure \ref{fig:paviaU} shows the scene for three different spectral reflectance bands 1, 50, and 103.





\section{Results}
\label{sec:results}

We have analyzed how different parameter settings affect the PSNR and computational times. 
We were able to show that different parameter values and their combination have an effect on the PSNR and computational times of the LRMR method.
The initial parameter values are optimized with respect to PSNR and CPU time.
The parameter value combinations are studied in detail and results are presented visually.

We study the initial parameter value combinations for three integer-valued parameters $r$, $s$, and $b$, and one real-valued parameter $p$.
We have applied nonlinear optimization with respect to the real-valued noise estimation parameter $p$ which resulted in an improvement in PSNR. 
The integer parameters have been analyzed on a series of test sets. The optimized values are chosen according to their best PSNR performance.
For all analyzed settings and resulting combinations, we show CPU times in Matlab and Python.
The computational experiments were performed on a laptop with Intel(R) Core(TM) i5-8250U CPU @ 1.60GHz 1.80 GHz processor and 16.0 GB RAM.

\begin{figure}[h]
    \centering
    \includegraphics[scale=0.9]{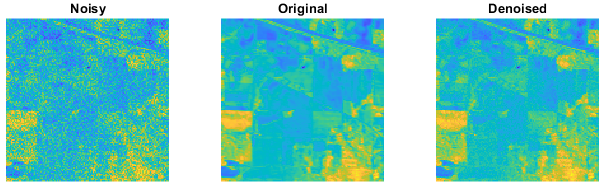}
     \caption{The noise removal performance of LRMR on a noise-corrupted spectral band for the Indian Pines data set. The image is restored efficiently.}
    \label{fig:lrmr}
\end{figure}

Figure \ref{fig:lrmr} shows the performance of the method for a noise-corrupted waveband of the Indian Pines data set \cite{IndianPines} for the 
parameter setting $r=7$, $p=0.15$, $b=20$, and $s=8$.

\subsection{Parameter analysis for $r$}

\begin{figure}[h] 
        \includegraphics[scale=0.4]{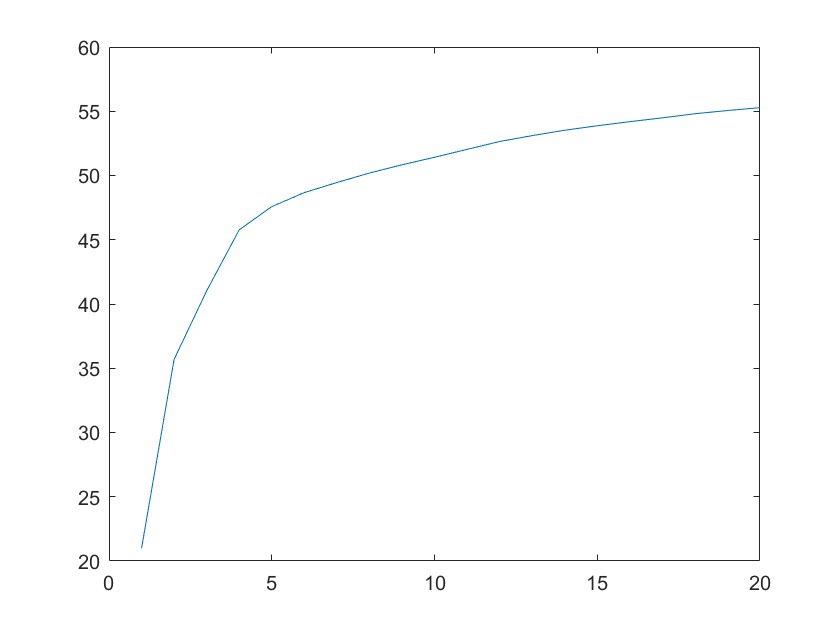}
        \caption{PSNR plot for different values of $r$.}
        \label{fig:psnr-r}
    \end{figure}

\begin{figure}[h] 
        \includegraphics[scale=0.4]{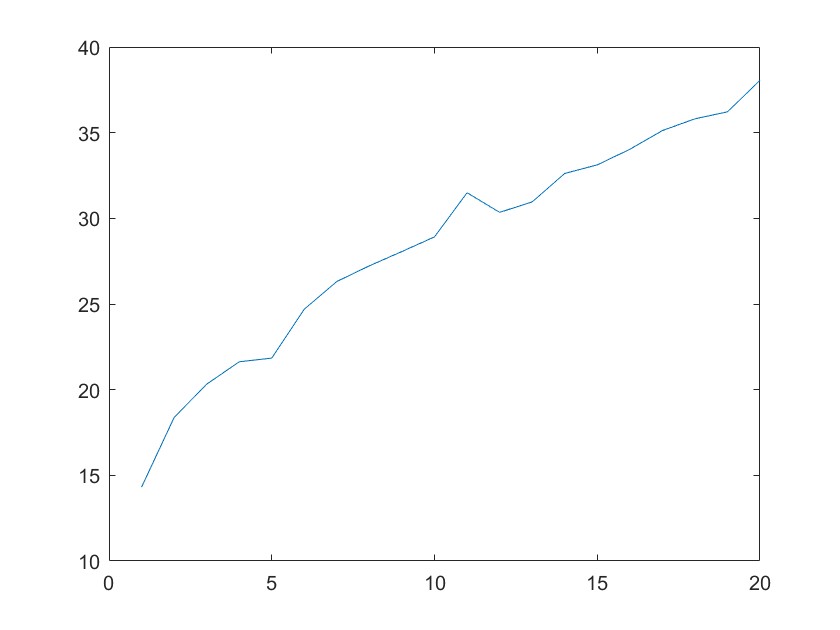}
       \caption{CPU time plot for different values of $r$.}
       \label{fig:time-r}
   \end{figure}

The rank parameter $r$ describes the upper bound for the rank of the low-rank matrix describing the noise-free signal of the datacube. 
The value should be as small as possible while not underestimating the noise intensity of the image data. 
Small values provide shorter computational times, while larger values provide higher PSNR values.
We vary the starting value for the rank parameter $r$ between 1 and 20. 
The other parameters have been kept $p=0.15$, $b=20$, and $s=8$.
Figure \ref{fig:psnr-r} shows the PSNR values for different values of $r$.
In general higher values for $r$ provide larger PSNR values. 
The analysis has shown that a value of $r \geq 4$ is sufficient. The largest gradient step is already taken from $r=1$ to $r=2$. 
The computational times in Matlab are shown in Figure \ref{fig:time-r}.
A larger rank value $r$ results in a larger computational time.
The computational tests were computed on the Indian Pines data set.
In the following, we present the PSNR, the gradient of the PSNR denoted as $\nabla$PSNR, and computational times in seconds in Matlab and Python denoted as $t_{Matlab}$ and $t_{Python}$, respectively, in Table \ref{tab:rank}.
Lower rank values are better in terms of computational effort.
\begin{table}[!ht]
\begin{center}
\begin{tabular}{|c|c|c|c|c|}
  \hline
   $r$ & PSNR & $\nabla$PSNR & $t_{\textit{Matlab}}$
   & $t_{\textit{Python}}$ \\
  \hline
   1  & 20.99 & 14.70& 14.31 & 95.15 \\
   2  & 35.69 & 10.02& 18.37 & 85.57 \\
   3  & 41.02 & 5.04 & 20.32 & 108.98 \\
   4  & 45.77 & 3.28 & 21.63 & 107.42 \\
   5  & 47.58 & 1.45 & 21.84 & 109.01 \\
   6  & 48.67 & 0.94 & 24.70 & 108.22 \\
   7  & 49.45 & 0.76 & 26.32 & 115.23 \\
   8  & 50.19 & 0.69 & 27.23 & 112.48 \\
   9  & 50.83 & 0.61 & 28.06 & 114.21 \\
   10 & 51.42 & 0.60 & 28.92 & 116.79 \\
   11 & 52.03 & 0.61 & 31.49 & 117.49 \\
   12 & 52.65 & 0.53 & 30.35 & 99.18 \\
   13 & 53.10 & 0.44 & 30.96 & 100.97 \\
   14 & 53.52 & 0.38 & 32.62 & 103.88 \\
   15 & 53.87 & 0.33 & 33.12 & 97.06 \\
   16 & 54.19 & 0.31 & 34.03 & 111.33 \\
   17 & 54.49 & 0.31 & 35.12 & 113.62 \\
   18 & 54.81 & 0.28 & 35.81 & 104.27 \\
   19 & 55.06 & 0.24 & 36.22 & 121.06 \\
   20 & 55.28 & 0.22 & 38.06 & 125.18 \\
  \hline
\end{tabular}
\end{center}
\caption{
The PSNR, gradient of PSNR, and the CPU times in seconds for computing the denoised cubes for different values of $r$ in Matlab $t_{\textit{Matlab}}$ and Python $t_{\textit{Python}}$.}
\label{tab:rank}
\end{table}
The computational effort is significantly more efficient in Matlab.
The analysis for Indian Pines suggests choosing $r=5$.
This parameter choice provides a sufficient balance between PSNR and computational time.

\subsection{Parameter analysis for $s$}

\begin{figure}[h] 
\includegraphics[scale=0.4]{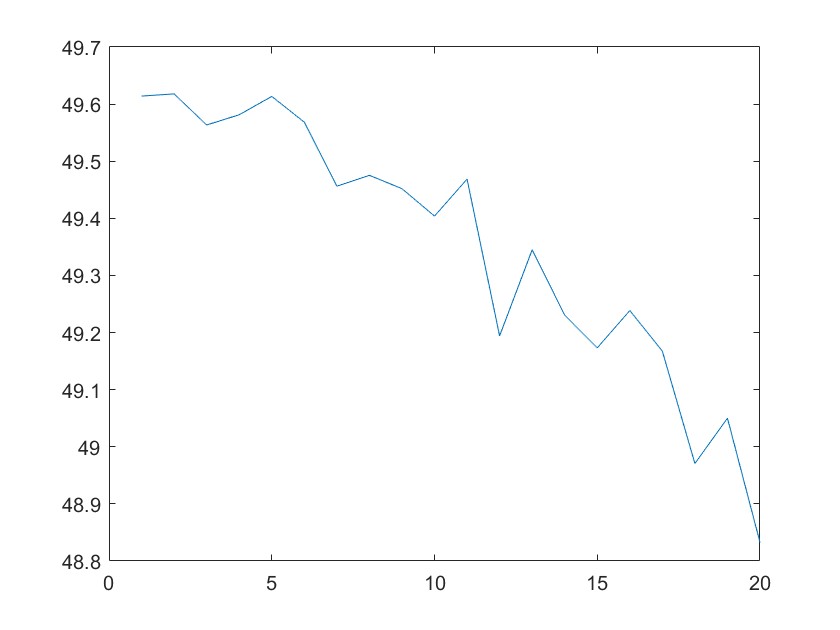}
        \caption{PSNR plot for different values of $s$.}
        \label{fig:psnr-s}
    \end{figure}

\begin{figure}[h] 
        \includegraphics[scale=0.4]{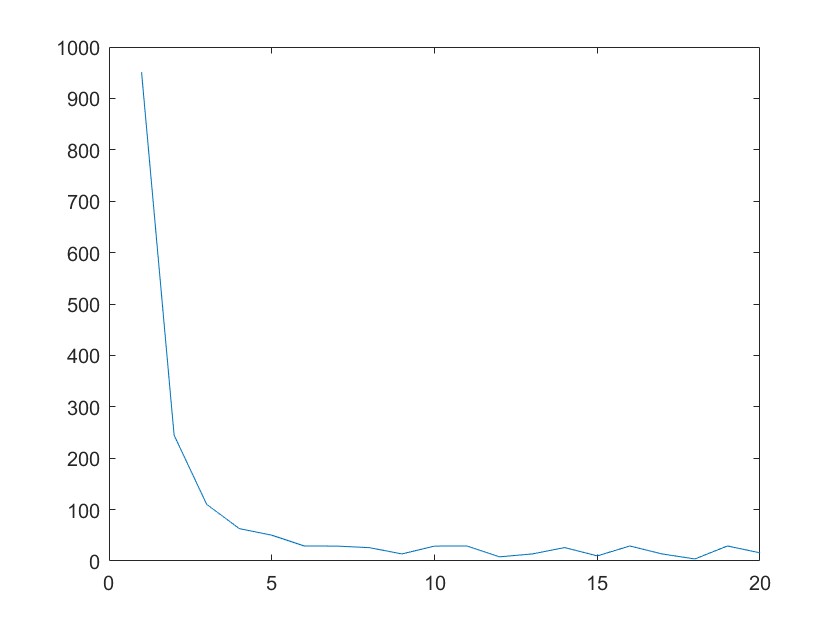}
        \caption{CPU time plot for different values of $s$.}
        \label{fig:time-s}
    \end{figure}

The parameter $s$ describes the iteration stepsize for processing the subcubes of LRMR's local image restoration.
We vary the parameter value for the stepsize $s$ between 1 and 20. The other parameter values have been kept $r=7$, $p=0.15$, and $b=20$.
Figure \ref{fig:psnr-s} shows the PSNR values for different values of $s$.
The trend shows that larger values for the stepsize $s$ yield smaller PSNR values. 
Hence, smaller values for $s$ are preferred, since we target a larger PSNR.
The computational times in Matlab are shown in Figure \ref{fig:time-s}.
The computational times are smaller for larger stepsizes $s$.
However, the CPU times decrease significantly for stepsizes larger than 4.   
The computational tests were computed on the Indian Pines data set.
In the following, we present the computational times in seconds in Matlab and Python in Table \ref{tab:stepsize}. 
Also, the gradient values for the computational times in Matlab $\nabla t_{Matlab}$ are presented.
The result is that a balance between smaller stepsizes that produce better PSNR and larger stepsizes that are better in terms of computational time and effort has to be achieved. 
However, note that the difference in PSNR values is not significant as Figure \ref{fig:psnr-s} shows. 
Figure \ref{fig:time-s} and the gradient values for the CPU times in Table \ref{tab:stepsize} show clearly that a stepsize $s \geq 4$ provides significant improvement regarding computational times.
\begin{table}[!ht]
\begin{center}
\begin{tabular}{|c|c|c|c|c|}
  \hline
   $s$ & PSNR  & $t_{\textit{Matlab}}$ & $\nabla t_{\textit{Matlab}}$
   & $t_{\textit{Python}}$ \\
  \hline
   1 & 49.61 & 951.11 & -706.38& 3793.27\\
   2 & 49.62 & 244.72 & -420.41& 1007.87\\
   3 & 49.56 & 110.29 & -90.82 & 481.01\\
   4 & 49.58 & 63.09  & -30.02 & 301.01\\
   5 & 49.61 & 50.24  & -16.93 & 193.88\\
   6 & 49.57 & 29.22  & -10.62 & 178.32\\
   7 & 49.46 & 29.00  & -1.63  & 114.31\\
   8 & 49.48 & 25.96  & -7.63  & 100.41\\
   9 & 49.45 & 13.74  &  1.54  & 139.66\\
   10& 49.40 & 29.04  &  7.71  & 112.07\\
   11& 49.47 & 29.16  & -10.48 & 111.66\\
   12& 49.19 & 8.08   & -7.71  & 137.42\\
   13& 49.34 & 13.74  &  9.01  & 53.73\\
   14& 49.23 & 26.11  & -1.90  & 100.47\\
   15& 49.17 & 9.94   &  1.53  & 38.31\\
   16& 49.24 & 29.18  &  1.85  & 111.84\\
   17& 49.17 & 13.63  & -12.60 & 52.35\\
   18& 48.97 & 3.98   &  7.81  & 178.98\\
   19& 49.05 & 29.25  &  5.95  & 114.47\\
   20& 48.83 & 15.88  & -13.36 & 60.94\\
  \hline
\end{tabular}
\end{center}
\caption{
The PSNR, the CPU times in seconds, and the gradient of CPU times for computing the denoised cubes for different values of $s$ in Matlab $t_{\textit{Matlab}}$ as well as the CPU times in Python $t_{\textit{Python}}$.}
\label{tab:stepsize}
\end{table}
As for the rank value $r$, the computational effort is significantly better in Matlab than in Python.

\subsection{Parameter analysis for $r$ and $s$}

\begin{figure}[h] 
\includegraphics[scale=0.4]{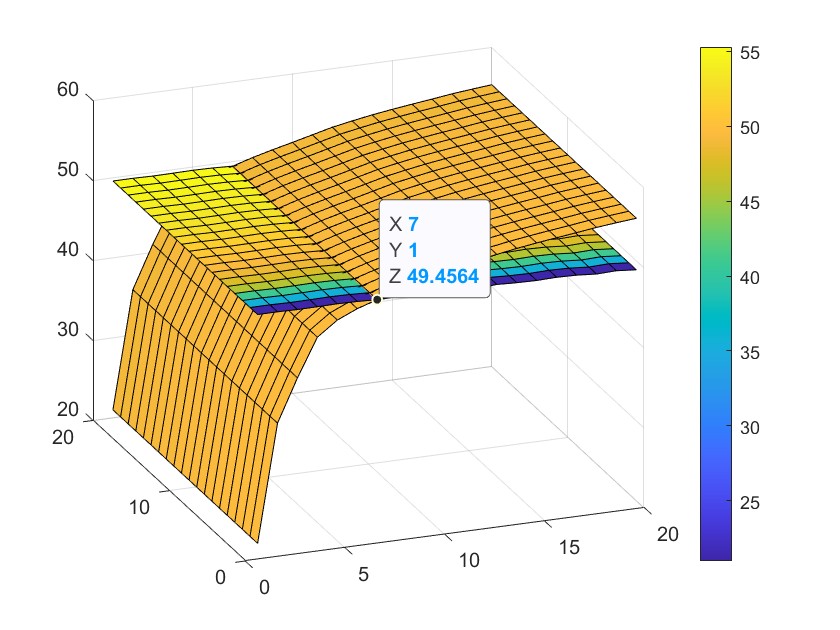}
        \caption{PSNR surface plots for different $r$ and $s$ intersect.}
        \label{fig:psnr-r-s}
    \end{figure}

In this parameter analysis, we study how the parameter choices for $r$ and $s$ influence each other.
Figure \ref{fig:psnr-r-s} shows the PSNR for different rank $r$ and stepsize $s$ values, whereas other parameter values $b$ and $p$ are set $b = 20$ and $p = 0.15$.
The point of intersection is marked and corresponds to $r = 7$ and values for $s \in [1,20]$.
Again $r=7$ has shown to be the optimized parameter choice, whereas values for $s$ show no significant difference in PSNR.
The tests were computed on the Indian Pines data set.

\subsection{Parameter analysis for $b$}

\begin{figure}[h] 
\includegraphics[scale=0.4]{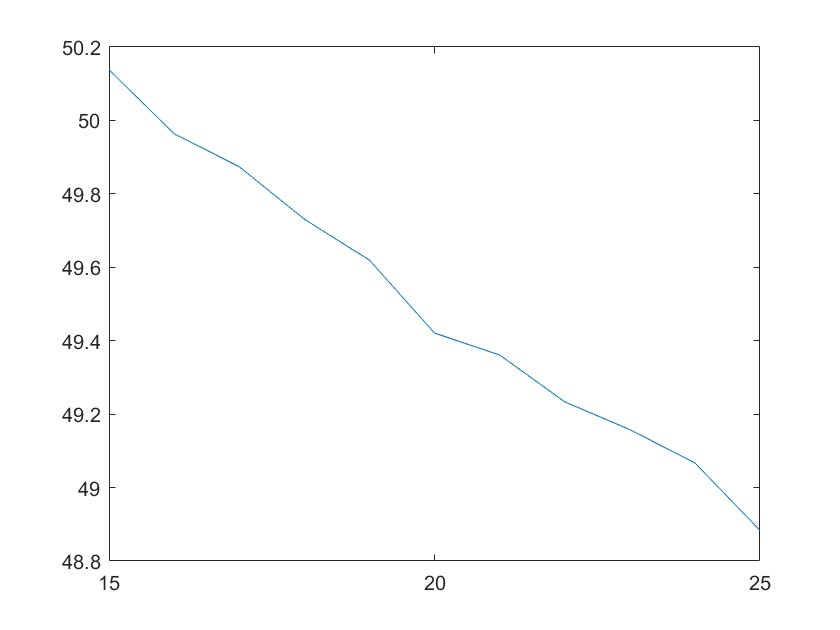}
        \caption{PSNR plot for different values of $b$.}
        \label{fig:psnr-b}
    \end{figure}

\begin{figure}[h] 
        \includegraphics[scale=0.4]{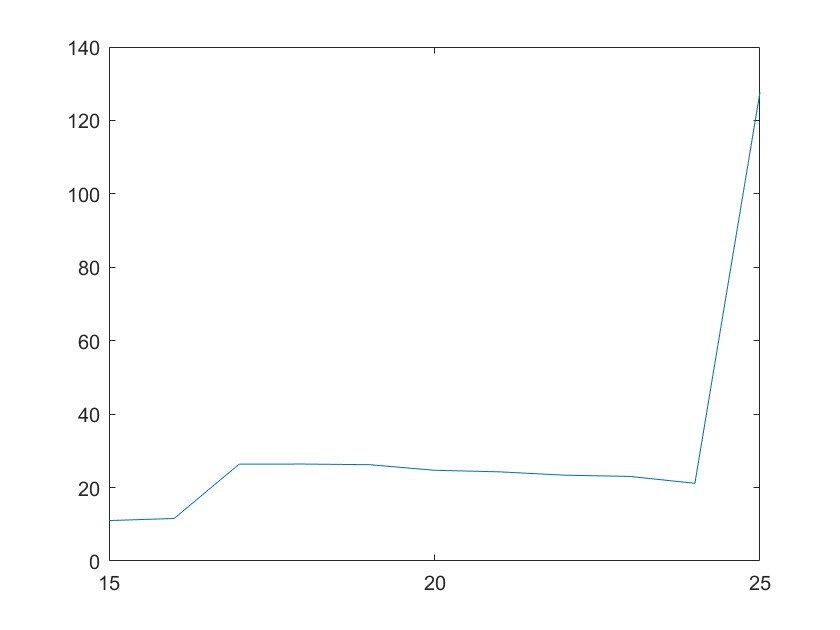}
        \caption{CPU time plot for different values of $b$.}
      \label{fig:time-b}
    \end{figure}

The blocksize of the subcubes for local restoration of LRMR is denoted by $b$.
The parameter choice is significant since naturally the size of the subcubes provides accuracy and computational effort to the restoration.
We vary the starting value for the blocksize parameter $b$ between 15 and 25.
Smaller values for $b$ have not provided a successful computation and error in PSNR.
The other parameters have been kept $r=7$, $p=0.15$, and $s=8$.
Figure \ref{fig:psnr-b} shows the PSNR values for different values of $b$.
The PSNR values are linearly decreasing with respect to increasing blocksize.
Figure \ref{fig:time-b} show the computational times in seconds in Matlab. Besides the large leap in the end and the smaller leap in the beginning, the computational time regarding the change in blocksize seems not to be affected significantly, although a weak trend towards increasing PSNR with increasing blocksize seems to be visibile.
We present the computational results 
for the Indian Pines data set.
In the following, we present the computational times in seconds in Matlab and Python in Table \ref{tab:block}.
\begin{table}[!ht]
\begin{center}
\begin{tabular}{|c|c|c|c|}
  \hline
   $b$ & PSNR  & $t_{\textit{Matlab}}$ 
   & $t_{\textit{Python}}$ \\
  \hline
   15 & 50.14 &  11.01 & 48.13\\
   16 & 49.96 &  11.58 & 83.56\\
   17 & 49.87 &  26.39 & 102.33\\
   18 & 49.73 &  26.41 & 124.62\\
   19 & 49.62 &  26.23 & 122.86\\
   20 & 49.42 &  24.71 & 127.33\\
   21 & 49.36 &  24.29 & 97.29\\
   22 & 49.23 &  23.37 & 92.63\\
   23 & 49.16 &  23.04 & 102.33\\
   24 & 49.07 &  21.16 & 179.84\\
   25 & 48.88 &  127.72& 257.69\\
  \hline
\end{tabular}
\end{center}
\caption{
The CPU times in seconds for computing the denoised cubes for different values of $b$ in Matlab $t_{\textit{Matlab}}$ and Python $t_{\textit{Python}}$.}
\label{tab:block}
\end{table}
Again the computational effort is better in Matlab.

\subsection{Parameter analysis for $p$}

\begin{figure}[h]
        \includegraphics[scale=0.4]{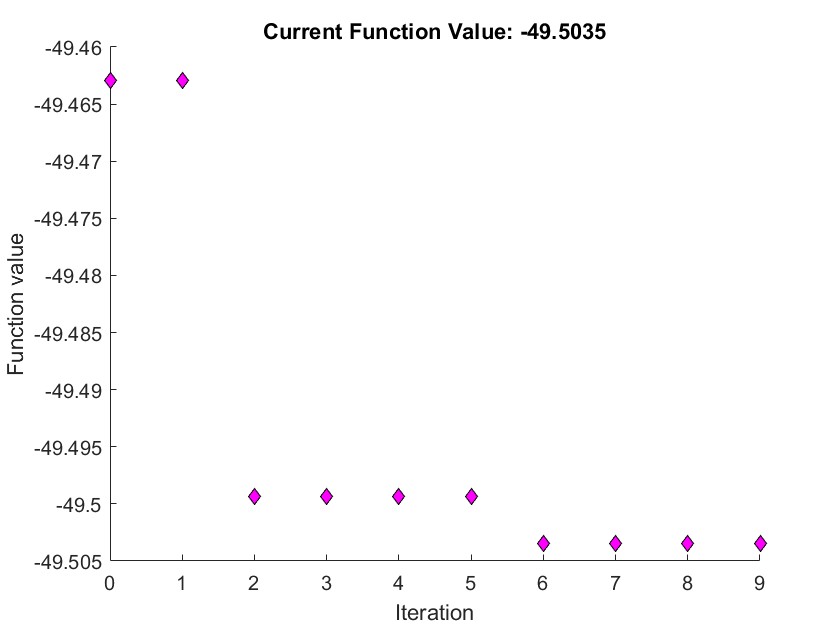}
        \caption{The function values of the iteration steps of the nonlinear optimization Nelder-Mead simplex algorithm with respect to $p$. For the Indian Pines data set and starting value $p=0.15$.} 
        \label{fig:optim}
    \end{figure}

The parameter $p$ describes an initial estimation of the amount of noise corruption in the datacube. 
It is the only real-valued parameter of LRMR's initial parameters. 
Hence, nonlinear optimization can be applied to determine the optimal initial value for $p$.
We present results for the 
nonlinear optimization method, the direct search method Nelder-Mead simplex algorithm \cite{Lagarias}.
We have applied the method on the different test sets and have varied the starting value for $p$. The other parameters were set $r=7$, $b=20$, and $s=8$.
We computed the optimization results of this section in Matlab.
In Figure \ref{fig:optim}, the negative PSNR values of the iteration steps of the nonlinear minimization method, the Nelder-Mead simplex algorithm, computed with respect to $p$ are presented. We present the results for the 
Indian Pines data set. 
The starting value $p = 0.15$ was chosen.
The results in Figure \ref{fig:optim} show convergence towards a local minimum.
The value yields an improvement in PSNR.
    The computational tests show a proper performance applying the Nelder-Mead simplex algorithm including convergence and computational effort. However, the method is searching for a local optimum and hence, has its limitations. The starting value is crucial. We include a table of final minimization function values and values for $p$ for different starting values, Table \ref{tab:p}.
    The stopping criteria where chosen as follows: the error in the function value is smaller than $1e-15$, the error in the optimization parameter value is smaller than $1e-9$, the maximum of function evaluations in one iteration step equals 20, and the overall maximum number of outer iteration steps equals 20.
\begin{table}[!ht]
\begin{center}
\begin{tabular}{|c|c|c|c|}
  \hline
   $p_0$ & -PSNR  & $p_{\textit{final}}$  \\
  \hline
   0.015& -49.990256 & 0.014977 \\
   0.05 & -49.887782 & 0.047813 \\
   0.15 & -49.503483 & 0.143438 \\
   0.5  & -49.476512 & 0.596875 \\
   1.   & -49.616882 & 0.975000 \\
   1.5  & -49.750781 & 1.818750 \\
  \hline
\end{tabular}
\end{center}
\caption{
The final function values for the minimization problem (the negative PSNR values) and the corresponding parameter value $p_{\textit{final}}$ for different starting values $p_0$. Clearly, local minima are approached. For the Indian Pines data set.}
\label{tab:p}
\end{table}   
We propose the choice of $p=0.15$ for further computations on the data sets used in this work. It is numerically stable and provides sufficient PSNR values.
Higher starting values than $1e+1$ have not provided sufficient convergence.

\begin{figure}[h] 
\includegraphics[scale=0.3]{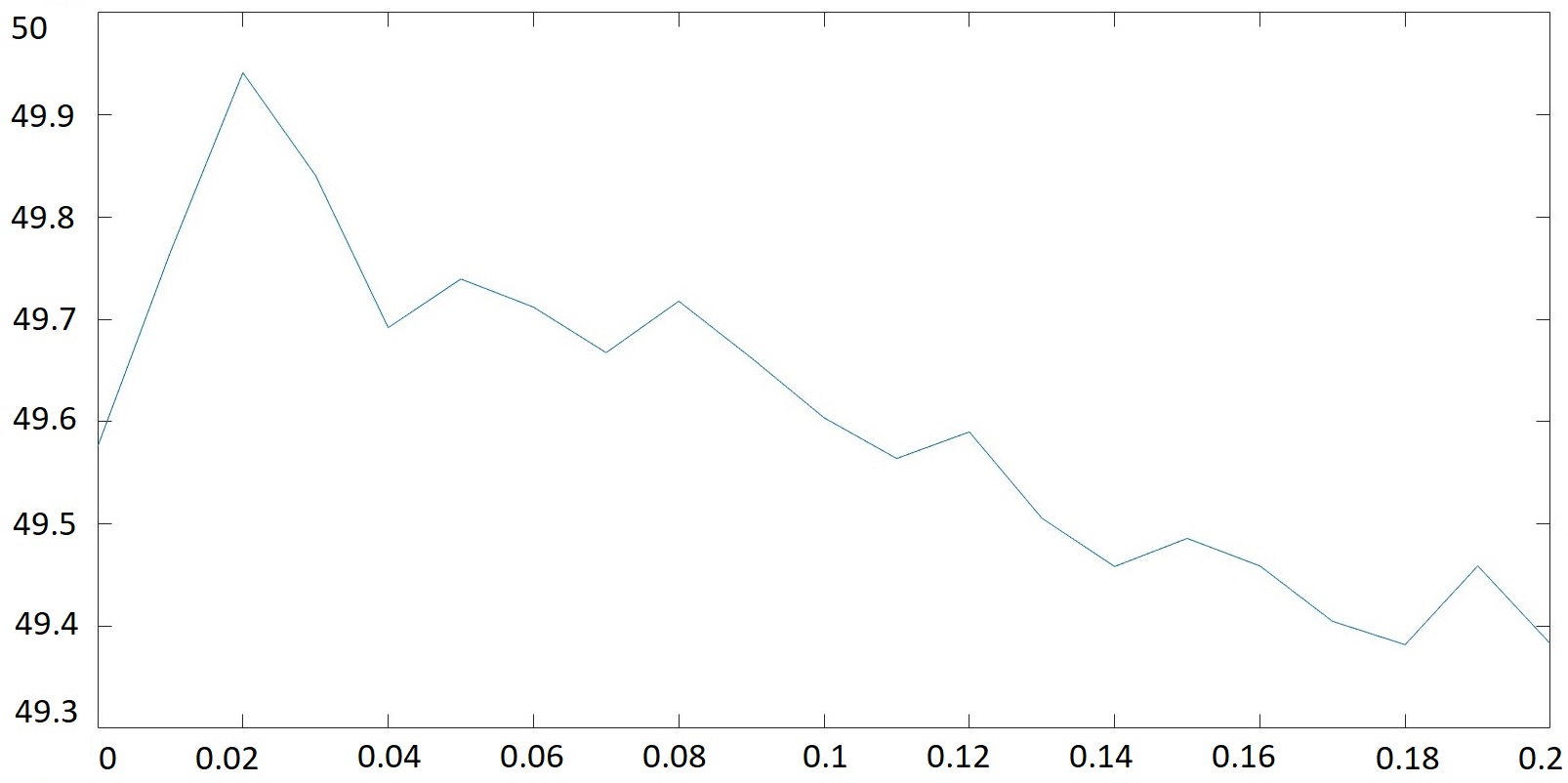}
        \caption{PSNR plot for analyzing different starting values of $p$.}
        \label{fig:psnr-p}
    \end{figure}

\begin{figure}[h] 
        \includegraphics[scale=0.4]{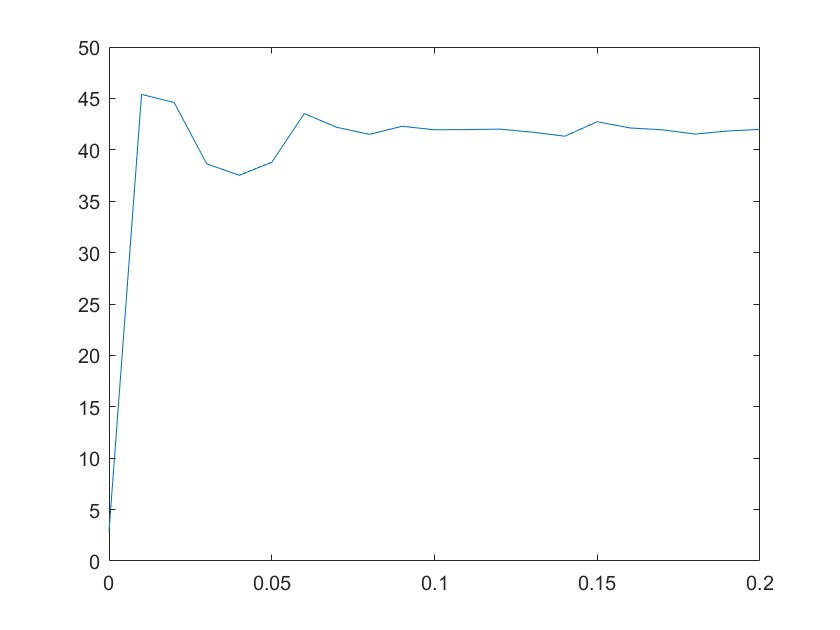}
        \caption{CPU time plot for different values of $p$.}
      \label{fig:time-p}
    \end{figure}

As second test, we have computed the PSNR and computational times for different values of $p \in [0, 0.2]$ to test different starting values.
The Indian Pines test set was used again for a qualitative parameter analysis.
Figure \ref{fig:psnr-p} shows the PSNR values for different $p$.
The figure suggests to choose a smaller starting value for $p$, which coincides with Table \ref{tab:p}.
Figure \ref{fig:time-p} shows the CPU times for different $p$.
The CPU times seem not be significantly affected in this range of $p$. 

\subsection{Optimized parameter choice}

Combining the results of the previous subsections, we present the resulting images for the optimized values.
We propose the following setting: $r=7$, $s=8$, $b=20$, $p = 0.15$.

\subsubsection{Indian Pines data set}

\begin{figure}[h]
    \centering
        \includegraphics[scale=0.35]{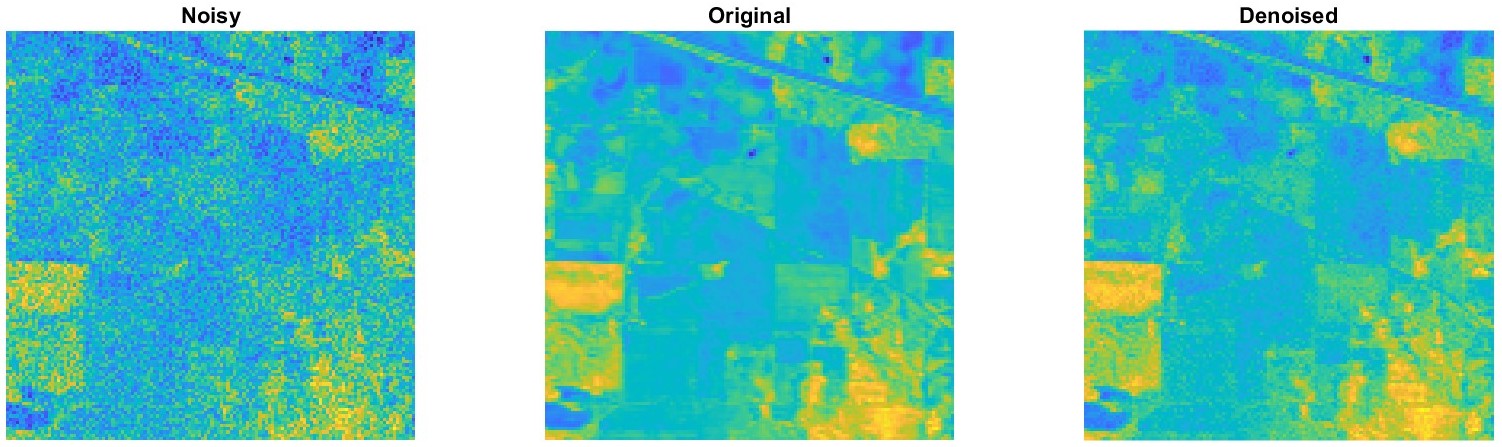}
        \caption{LRMR result computed on the Indian Pines data set for the optimized parameter setting.}
        \label{fig:indian_final}
  \end{figure}

Figure \ref{fig:indian_final}.

\subsubsection{Pavia Centre data set}

\begin{figure}[h]
    \centering
    \includegraphics[scale=0.35]{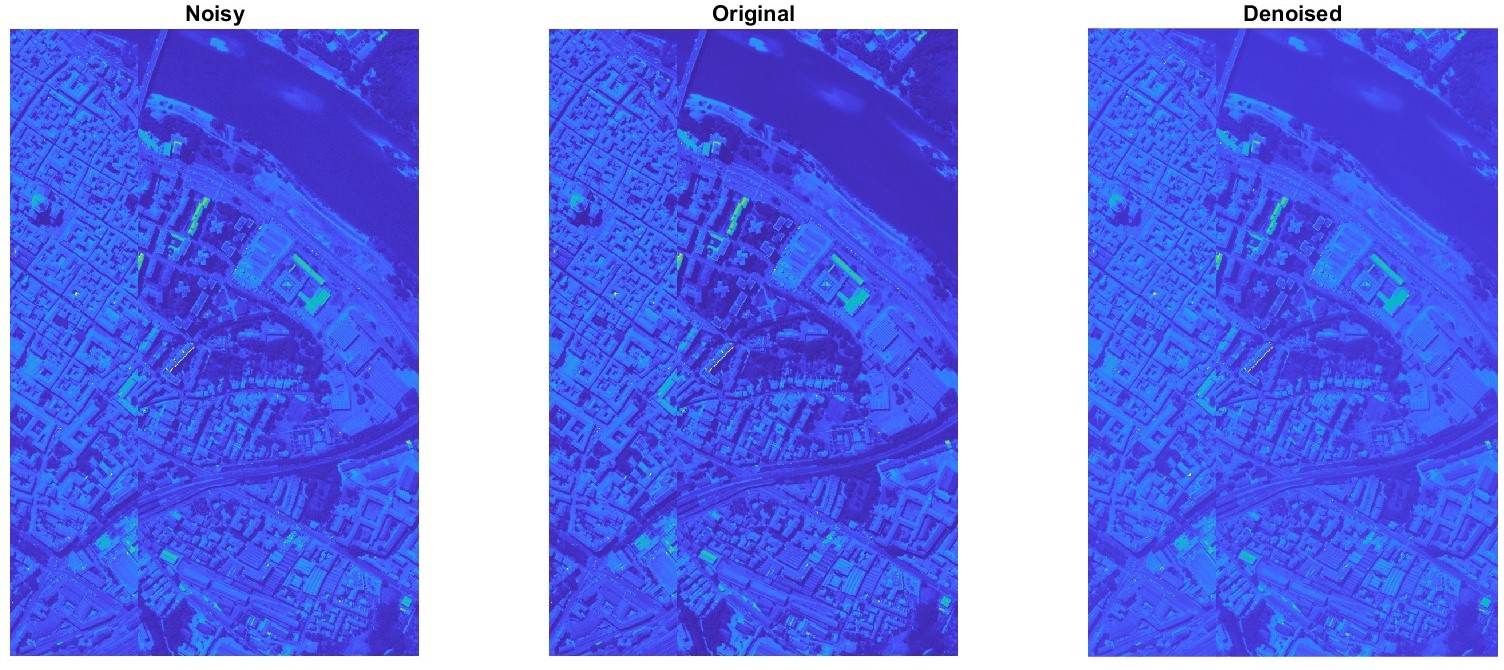}
        \caption{LRMR result computed on the Pavia Centre data set for the optimized parameter setting.}
        \label{fig:pavia:final}
\end{figure}

Figure \ref{fig:pavia:final}.

\subsubsection{Pavia University data set}

\begin{figure}[h]
    \centering
    \includegraphics[scale=0.35]{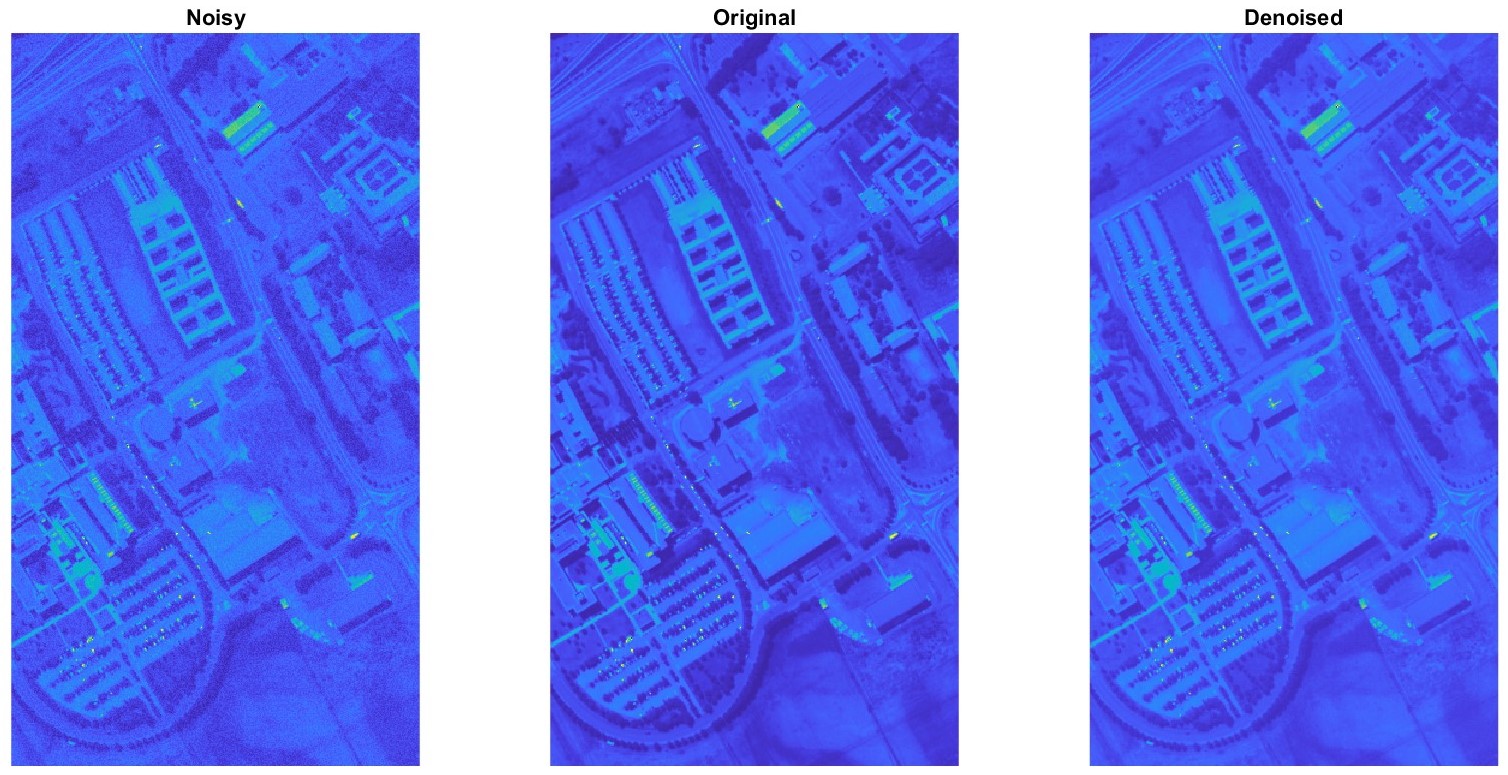}
        \caption{LRMR result computed on the Pavia University data set for the optimized parameter setting.}
        \label{fig:paviaU:final}
\end{figure}

Figure \ref{fig:paviaU:final}.

\section{Discussion}
\label{sec:discussion}

We note that parameter tests on other different hyperspectral image data sets may provide a different parameter setting. Hence, we advise to perform a parameter study for the initial values in the LRMR method set up for a new test set.
A problem specific parameter configuration is suggested and it was a goal of this work to highlight and provide a sample analysis as model for other tests.

We computed the nonlinear optimization and the main analysis results in Matlab.
A comparison between the computational costs of Matlab and Python has been presented for the integer valued parameters. The overall behavior of the algorithms is the same in Matlab and Python. However, the computational times have been significantly larger in Python.
The reason is the efficient storing and processing of matrix-vector operations in Matlab.
The analysis of computational times is important for applications, where the data budget is limited or where the data has to be processed in real-time, therefore for real-time decision-making.
It was important to show the similar behavior in relative terms of both implementations, Matlab and Python, and to discuss the implementation issues that may occur. We observed a stable behavior for each implementation and the relative behavior was the same.

Regarding future work, other optimization methods could be applied in order to select the optimized initial values for $p$. We have applied the Broyden–Fletcher–Goldfarb–Shanno (BFGS) Quasi-Newton method \cite{Broyden, Fletcher, Goldfarb, Shanno} as alternative to the nonlinear optimization. However, the method stopped always after two iteration steps. Alternatively, the use of a machine learning approach might provide additional parameter optimization results and a more efficient optimization approach to choosing $p$. 
This would require access to or the creation of a large data set in order to secure efficient training and testing of the machine learning algorithm.
However, this is beyond the scope of this article and is ongoing work of the author.


\section{Conclusions}
\label{sec:conclusions}

The parameter analysis and optimization yield improvement of the hyperspectral images regarding the PSNR.
We present 
a detailed analysis with respect to the initial parameters of the LRMR method for different data test sets, however, with a focus on the Indian Pines data set. 
The approach of this work studying the parameter optimization for LRMR's initial parameters is new and has not been presented in that way earlier.
The analysis provides new insights into the flexibility and boundaries of the method. 


\section*{Acknowledgment}

This research was funded by the Regional Council of Central Finland/Council of Tampere Region and European Regional Development Fund as part of the \textit{coADDVA - ADDing VAlue by Computing in Manufacturing} projects of Jamk University of Applied Sciences.




\bibliographystyle{abbrv}
\bibliography{bibliographyAS}

\end{document}